\documentclass[smallcondensed, final]{svjour3}

\usepackage{t1enc}

\usepackage{amsmath,amssymb}
\usepackage{color}
\usepackage{ifpdf}

\ifpdf
\usepackage[pdftex]{graphicx}
\else
\usepackage{graphicx}
\fi

\newcommand{\dd}{\; \mathrm{d}}
\newcommand{\bs}[1]{\boldsymbol{#1}}

\newcommand{\N}{\mathbb{N}}
\newcommand{\R}{\mathbb{R}}

\newcommand{\Z}{\mathbb{Z}}
\newcommand{\ul}{\underline}
\newcommand{\ol}{\overline}
\newcommand{\rwwli}{RWwLI }

\renewcommand{\l}{\left}
\renewcommand{\r}{\right}

\spnewtheorem*{proof*}{}{\it}{\rm}

\begin{document}
\title{Locally Perturbed Random Walks with Unbounded Jumps}
\author{Daniel Paulin \and Domokos Sz\'asz}
\institute{Daniel Paulin \at Budapest University of Technology and Economics, Hungary and National University of Singapore, Singapore\\
\email{paulindani@gmail.com}
\and Domokos Sz\'asz \at Department of Stochastics, Budapest University of Technology and Economics, Hungary\\\email{szasz@math.bme.hu}
}
\begin{abstract}
In \cite{SzT}, D. Sz\'asz and A. Telcs have shown that for the diffusively scaled, simple symmetric random walk, weak convergence to the Brownian motion holds even in the case of local impurities if $d \ge 2$. The extension of their result to finite range random walks is straightforward. Here, however, we are interested in the situation when the random walk has unbounded range. Concretely we generalize the statement of  \cite{SzT} to unbounded random walks whose jump distribution belongs to the domain of attraction of the normal law. We do this first: for diffusively scaled random walks  on $\mathbf Z^d$  $(d \ge 2)$ having finite variance; and second: for random walks with distribution belonging to the non-normal domain of attraction of the normal law.  This result can be applied
to random walks with tail behavior analogous to that of the infinite horizon Lorentz-process; these, in particular, have infinite variance, and convergence to Brownian motion holds with the superdiffusive $\sqrt{n \log n}$ scaling.
\end{abstract}
\maketitle

\section{Introduction}

Our goal in this paper is to show that local impurities do not influence the - appropriately scaled - weak limit behavior of random walks on $\mathbf Z^d$ $(d \ge 2)$.
In \cite{SzT}, D. Sz\'asz and A. Telcs have shown that for the diffusively scaled, simple symmetric random walk, weak convergence to the Brownian motion holds even in the case of local impurities if $d \ge 2$. The extension of their result to finite range random walks is straightforward. Here, however, we are interested in the situation when the random walk has unbounded range (for simplicity, we always assume that the distribution of the jumps is centered, i.e. their expected value is 0). Concretely we generalize the statement of  \cite{SzT} to unbounded random walks whose jump distribution belongs to the domain of attraction of the normal law. We do this first: for diffusively scaled random walks  on $\mathbf Z^d$ $d \ge 2$ having finite variance, and second: for random walks with distribution belonging to the non-normal domain of attraction of the normal law (we note that, if $d=1$, this domain consists of distributions which have infinite variance but $L(x)= \int_{|u| \le x} u^2 dF(u)\to \infty $ is a slowly varying function; in this case partial sums converge to the gaussian law by using the slightly superdiffusive scaling $B_n$, where $B_n$ is determined by the relation $\frac{nL(B_n)}{B_n^2}\rightarrow 1$. For the reader's convenience, a summary of related definitions and results can be found in the Appendix). 

The idea of \cite{SzT} served as an intuitive background for proving that the diffusively scaled, locally perturbed, planar, finite horizon Lorentz process
converges weakly to the Brownian motion (cf. \cite{ChD} and \cite{DSzV}). Our second result is hoped to provide a starting point to prove the corresponding statement for the infinite horizon Lorentz process, recently a widely studied object (cf. \cite{SzV2007}, \cite{MS08}, \cite{ChD2009},  \cite{ChD2010}). (In fact, for stochastic models of the infinite horizon Lorentz processes other types of perturbations should also be taken into account, an object of future research.)

In Section 2, we are going to state our general theorems, Section 3 contains some lemmas and definitions, and Section 4 proves the theorems using these. The proofs of the lemmas can be found in Section 5. Finally, Section 6 is devoted to comments.

\section{Main results (in particular, the generalization of  \cite{SzT} to the case of infinite horizon)}
The main difficulty in generalizing Theorem 1 of \cite{SzT} to the case of random walks with unbounded jumps is that in this case the coupling argument of \cite{SzT} breaks down. Though the structure of our proofs is similar to that used in \cite{SzT}, for avoiding the aforementioned difficulty one needs novel ideas (as seen, in particular, in the proofs of our lemmas in Section 5).

To avoid unnecessary complications we suppose throughout the whole paper: (i) the dimension $d \ge 2$; (ii) the aperiodicity of the random walk (i. e. if we denote by $P_n(z \to x)$ the $n$-step transition probability of the random walk, then we assume: if there is an $m$ such that $P_m(z\to x)>0$, then we can choose $n_0$ (which may depend on $x$ and $z$) such as for all $n>n_0$, $P_n(z \to x)>0$). Finally  $|| . ||$ denotes  the maximum ($L_\infty$) norm for vectors.

\begin{definition}\label{D:RW}
Let  $y_1, y_2, \ldots $ be independent, identically distributed (or briefly iid) random variables with
\[P(y_i= y)=\ul{P}(y) \hspace{6mm} y\in \Z^d.\]
The stochastic process $Y_0, Y_1, \ldots$ defined by
\[Y_n=y_0+\sum_{i=1}^n y_i\]
where $y_0=z \in \Z^d$, is a \emph{random walk} (or briefly RW). The measure defined by this random walk will be denoted by $\ul{P}_z$, and the transition probabilities of the process will be denoted by  $\ul{P}(x,y)$ ($\ul{P}(x,y)=\ul{P}(0,y-x)\overset{def}{=}\ul{P}(y-x)$).
\end{definition}

As usual, the random walk orbit determines a continuous time stochastic process $\eta^Y_n: [0,1]\to \R$ with continuous trajectories in the following way: $\eta^Y_n(t)=n^{-1/2}Y_{nt}$ if $t=0,\frac{1}{n},\frac{2}{n}, \ldots, 1, \ldots$, and it is linear between these points. It is well known that if $\ul{P}$ is such that second moments of the jumps are finite, i.e.
\[\underset{x\in \Z^d}{\sum}\ul{P}(x)||x||^2<+\infty,\] then, as $n\to \infty$, $\eta^Y_n$ converges weakly to $W_{\sum}(t)$, a $d$-dimensional Wiener-process in the space $C_d[0, \infty)$ (with $d\times d$ covariance matrix $\sum$ equal to that of the random vector $y_i$). Weak convergence in $C_d[0, \infty)$ will be denoted by $\Rightarrow$. The proof of this result ($d$-dimensional Donsker functional central limit theorem) can be found in \cite{Whitt} page 106, Theorem 4.3.5.
\begin{remark}
Since our arguments that show weak convergence in $C_d[0, 1]$ also imply that in $C_d[0, \infty)$ in a  standard way (cf. \cite{Lindvall}), we often will only formulate our statements for $C_d[0, 1]$.
\end{remark}

\begin{definition} \label{D:ISRW2}
Let P be a transition probability matrix on $\Z^d$ (so $P(x,y)$ is not necessarily equal to $P(0,y-x)$). We call the time-homogeneous Markov process $X_i$, $i=0,1\ldots$ with transition probabilities $P$ a \emph{random walk in an inhomogeneous medium}.
\end{definition}

\begin{definition}\label{D:RWwLI2}
If $X_n$ is a random walk in an inhomogeneous medium and there exists a finite set $A\subset \Z^d$ such that for all $u\notin A$, $v\in \Z^d$, $P(u,v)=\ul{P}(u,v)$, then we call $X_n$ a \emph{random walk with local impurities} (or briefly \rwwli).
\end{definition}

Let us define a directed graph $G=(\Z^d, E)$, where $(E=\{(u,v)| P(u,v)\ne 0\})$.
In both our theorems we assume:

\noindent {\bf Assumptions:}
\begin{enumerate}
\item [(i)]
let $X_n$ be a \rwwli, where the starting point $z$ lies in the infinite, strongly connected component $Q$ of $G$ (a directed graph $Q$ is  strongly connected  if there is a directed path from each vertex to each other vertex in $Q$, a strongly connected component of $G$ is a maximal strongly connected subgraph); for simplicity, assume also that $0 \in Q$;
\item [(ii)]
assume that there is an $\epsilon>0$ that for all the impurities, the jump from the impurity has a distribution whose $\epsilon$th moment exists (i.e. $E(||J||^{\epsilon})<+\infty$, where $J$ is the jump from an impurity). 
\end{enumerate}

\begin{theorem}\label{the1}
Define $\eta^X_n: [0,1]\to \R$ analogously to $\eta^Y_n(t)$. 
Suppose that the second moment of the jumps (with respect to $\ul{P}$) exists, their expected value is 0 and the RW is aperiodic.
Then, as $n\to \infty$,
\[\eta^X_n(t)\Rightarrow W_{\sum}(t), \hspace{6mm} t\in [0,1].\]
\end{theorem}

\begin{definition}\label{D:Btype}
We call a random variable (or random vector) $X$ with $\mathbb E||X||^2 = \infty$ B-type if its expected value is 0, and it belongs to the domain of attraction of the normal law. (As a consequence, for partial sums of iid B-type summands the scaling is larger than $\sqrt{n}$, i.e. $\underset{n\to \infty}{\lim} \frac{B_n}{\sqrt{n}}=\infty$). The reader is reminded that some facts from classical limit theory of probability, e. g. the definition of domain of attraction,  are collected in the Appendix.
\end{definition}

\begin{theorem}\label{the2}
Suppose that the distribution of the RW jumps is B-type with scaling $B_n$, their expected value is 0 and the RW is aperiodic. Define  $\eta^{X'}_n: [0,1]\to \R$ as
\[\eta^{X'}_n(t)=B_n^{-1}X_{nt} \text{ if } t=0,\frac{1}{n},\frac{2}{n}, \ldots, 1,\] and let it be linear between these points.
Then as  $n\to \infty$,
\[\eta^{X'}_n(t)\Rightarrow W_{\sum}(t), \hspace{6mm} t\in [0,1].\]
\end{theorem}

\begin{remark} 
Of course, it may also occur that partial sums of iid random vectors would 'naturally' scale differently in different directions. This is, for instance, the case in Theorem 7 of \cite{ChD2009}: they consider a two-dimensional periodic Lorentz process where all the collision-free trajectories are parallel to the $x$-axis and the scaling is $\sqrt{n \log n}$ in the direction of the $x$-axis whereas it is $\sqrt{n }$ in the direction of the $y$-axis. The reader can convince him/herself that general case is analogous: in some directions one has to use the strongest scaling and then in some orthogonal directions the next to the strongest ones, etc. and our methods are also applicable in this situation.
\end{remark}

\begin{remark}
For $d=1$, Theorem 1 is not true; cf. \cite{Skew}.
\end{remark}
\begin{remark}
In Theorem \ref{the1}, $\sum$ is the same as the covariance matrix of the jumps of the RW (having distribution $\ul{P}$). For Theorem 2, $\sum$ is still determined by the distribution of the RW jumps (see Theorem 4.2 of \cite{Rv} in the Appendix).
\end{remark}

\section{Preliminary notes to the proofs of Theorems \ref{the1} and \ref{the2}}
Before the rigorous discussion we want to show the idea of the result. If $d\ge 3$, then for the simple symmetric random walk, P\'{o}lya's theorem says that with a probability 1, the number of returns into the origin --- or into a finite set $A$ --- is finite. This can be simply proven for any non-degenerate random walk. For instance, for $d=3$, once the Local Central Limit Theorem (\cite{Lawler}, page 25) holds, stating that the probability of return into the origin in the $n$th step is $p_n=O(n^{-3/2})$, then the series $p_i$ is summable. Thus the expected number of returns is finite. Therefore, in case of a finite modification, the random walk leaves the set $A$ after a finite time. Consequently, in the limit, the effect of the modification vanishes.

Much more interesting is the case $d=2$, when the random walks we are interested in are recurrent. Again, once for some sequence $\{B_n > 0\}_{n \ge 1}$ \ $\frac{Y_n}{B_n}$ has a limit law, then under some additional conditions, its local version also holds. It implies that the expected number of returns into a finite set until time $n$ is $O(\sum_{j=0}^n \frac{1}{B_j^2})=O(\sum_{j=0}^n \frac{1}{j})$ which is always $O(\log n) = o(\sqrt n)$.  The expected time spent in $A$ during one visit is uniformly bounded and since the normalizing factor is at least of order $\sqrt{n}$, the previous conclusion is also true.

Now, we are going to start with definitions and our key lemma (its proof can be found in section 5).

\begin{definition}
The impurities are local, so we can select and fix an $N\in \N$ such that all the impurities lie in the cube \[K_N=\l[-N-\frac{1}{2}; N+\frac{1}{2}\r]^d \subset \R^d \hspace{5mm}\text{i.e., } A\subset K_N.\]
\end{definition}

\begin{definition}
\[\rho_n=\sum_{i=0}^n \bs{1}_{X_i\in K_N}\]
i.e. the time spent by the \rwwli in $K_N$ until time $n$.
\end{definition}

\begin{lemma}\label{rhoon12}
If $d\ge 2$, and $z\in Q$, then, as $n\to \infty$, \[E_z(\rho_n)=O(\log n).\]
\end{lemma}

\section{Proofs of Theorems \ref{the1} and \ref{the2} assuming Lemma \ref{rhoon12}}
\begin{proof*}{\bf Proof of Theorem \ref{the1}}
We are going to write the proof for dimension $d\ge 2$ (even though the $d\ge 3$ case could be proven simply, as indicated before). With the help of $X_n$ we define a new process $Z_n$. Let $Z_0=z$ and
$Z_{n+1}-Z_n=X_{n+1}-X_n$ if $X_n\notin K_N$, while if $X_n\in K_N$, then let $Z_{n+1}-Z_n$ be independent of
$X_0, X_1, \ldots X_n$ and of $Z_0, Z_1, \ldots Z_n$, and let $P(Z_n,Z_{n+1})=\ul{P}(Z_n,Z_{n+1})$. Then it is clear that $Z_n$ is a random walk with transition matrix $\ul{P}$. Let's define $\eta^Z_n(t)$ analogously as $\eta^Y_n(t)$ in Section 2, then \[\eta^Z_n(t)\Rightarrow W_{\sum}(t)\]
Thus in order to establish that \[\eta^X_n(t)\Rightarrow W_{\sum}(t)\] also holds, because of the piecewise linearity of $\eta^X_n(t)$ and $\eta^Z_n(t)$, it is sufficient to show that \[n^{-1/2}\cdot \underset{t\in[0,1]}{\sup} ||X_{[nt]}-Z_{[nt]}|| \to 0 \text { as } n\to \infty.\]

Observe that
\begin{align*}
n^{-1/2}\cdot \underset{t\in[0,1]}{\sup} ||X_{[nt]}-Z_{[nt]}||&\le n^{-1/2}\cdot \sum_{i=1}^n ||(X_i-X_{i-1})-(Z_i-Z_{i-1})||=\\
&=n^{-1/2} \cdot \sum_{i: X_{i-1}\in K_N} ||(X_i-X_{i-1})-(Z_i-Z_{i-1})||
\end{align*}

Let us define a sequence of random variables $J_k=||k \text{th step from }K_N||$, $k\in \Z^+$. Then our task is to show that
\[n^{-1/2} \cdot \sum_{k=1}^{\rho_n} J_k\Rightarrow 0\]
Firstly, we will show this for the case when $K_N$ consists of only one point, and then we show it for the general case.

In the one - point case, $J_k$ are independent and they have the same distribution, which, by assumption, has $E(J_k^{\epsilon})<\infty$ for some $\epsilon>0$, by a simple argument we can also suppose that $\epsilon<1$. Let us call $E(J_k^{\epsilon})=K$, then by Markov - inequality, we have $P(J_k^{\epsilon}>h)<K/h$ ($h>0$), therefore $P(J_k>h)<K/h^{\epsilon}$. Let us define a new random variable $J'$ as $P(J'>h)=1$ if $h\le 0$ and $P(J'>h)=\min(1, K/h^{\epsilon})$ if $h>0$. Then we define the sequence $J_k'$ with the same distribution as $J'$, and can see that for every $\gamma \ge 0, $\[P\l(n^{-1/2} \cdot \sum_{k=1}^{\rho_n} J_k > \gamma\r) \le P\l(n^{-1/2} \cdot \sum_{k=1}^{\rho_n} J_k' > \gamma\r)\]
therefore it is sufficient to show that the latter converges to zero as $n \to \infty$.
Then we can easily verify that this distribution of $J'$ belongs to the domain of attraction of stable law with parameter $\alpha=\epsilon$. Therefore if we sum $k$ random variables of this distribution, and note the sum by $S_k$, then $\frac{S_k}{k^{1/\epsilon}} \Rightarrow S(\epsilon, \beta)$, where $S(\epsilon, \beta)$ is a stable law.

Now we can write
\begin{align*}
P\l(n^{-1/2} \cdot \sum_{k=1}^{\rho_n} J_k' > \gamma\r)&\le P\l(\l.n^{-1/2} \cdot \sum_{k=1}^{\rho_n} J_k' > \gamma\r|\rho_n>\log^2(n)\r)\cdot P\l(\rho_n>\log^2(n)\r)+\\ &+P\l(n^{-1/2} \cdot \sum_{k=1}^{[\log^2 n]} J_k' > \gamma\r)
\end{align*}
By Lemma \ref{rhoon12}, $E(\rho_n)=O(\log n)=o(\log^2 n)$. Hence the first term tends to zero by the Markov - inequality. Therefore we only need to show that the second term tends to zero too.
\[P\l(n^{-1/2} \cdot \sum_{k=1}^{[\log^2 n]} J_k'>\gamma\r)=P\l(\l(\log^{-2/{\epsilon}}(n) \cdot \sum_{k=1}^{[\log^2 n]}J_k\r) \cdot \frac{n^{-1/2}}{\log^{-2/{\epsilon}}(n)}>\gamma\r)\]
We have $P\l(\log^{-2/{\epsilon}}(n) \cdot \sum_{k=1}^{[\log^2 n]}J_k>x\r) \to P\l(S(\epsilon,\beta)>x\r)$ for every $x>0$ as $n\to \infty$, so it is obvious that $P\l(n^{-1/2} \cdot \sum_{k=1}^{[\log^2 n]} J_k' > \gamma\r)\to 0$ as $n\to \infty$.

The second case is when $K_N$ consists of several points. First  we denote by $V_n$ the set ${x_1,x_2, \ldots x_n}$ of the first $n$ points in $K_N$ hit by $X_j$. Let $W_n$ be the set of all possible $V_n$ that occur with probability greater that zero. Then we can write \[P\l(n^{-1/2} \cdot \sum_{k=1}^{\rho_n} J_k > \gamma\r)=
\sum_{V_n\in W_n} P\l(\l.n^{-1/2} \cdot \sum_{k=1}^{\rho_n} J_k > \gamma \r| V_n\r) P(V_n)\]
Now it is sufficient for us to show that $P\l(\l.n^{-1/2} \cdot \sum_{k=1}^{\rho_n} J_k > \gamma \r| V_n\r) \to 0$ as $n\to \infty$ for every $V_n\in W_n$.

If we fix $V_n$, then $J_k$ are conditionally independent. Moreover, we can see that $J_k|V_n = J_k|\{x_k,x_{k+1}\}$, so it only really depends on the $k$th and $k+1$th step in $K_N$. For fixed $x_k$,
$J_k|x_k$ has finite $\epsilon$th moment for some $\epsilon>0$ (independent of $k$ and $x_k$). If \[K=\underset{x\in K_N \cap Q}{\max} E(J_k^{\epsilon}|x_k=x),\] then we see that $P\l(J_k>h|x_k\r)<K/h^{\epsilon}$.
But \[P(J_k>h|x_k)=\sum_{x\in K_N} P(x_{k+1}=x|x_k)\cdot P(J_k>h|x_k,x_{k+1}).\]
$P(x_{k+1}=x|x_k)$ is independent of $k$, it only depends on the two points in $K_N$, so it has at most $|K_N|^2$ different values. We only need to consider cases when it is not zero, because if it is zero, than it can never occur in any $V_n\in W_n$. Therefore there is a positive $p>0$ such that $P(x_{k+1}=x|x_k)>p$ if it is not zero. So
\[P(J_k>h|V_n)<\frac{K}{ph^{\epsilon}}\]
From here, because of the conditional independence of $J_k$, we can follow the proof of the one - point case.
\end{proof*}

\begin{remark}
The condition of existence of $\epsilon$th moment seems difficult to avoid, since if we take the case of a distribution for $J$ such that the cumulative distribution function has the form $F(x)= 1-\frac{1}{\log^2(e+x)}$ for $x>0$ and 0 for $x\le 0$, then we can observe the sum of $\log n$ independent random variables of this kind, divided by $\sqrt{n}$, will not converge to zero. Nevertheless, it is possible that in this case the random walk returns less frequently than $O(\log n)$ times in $n$ steps, but this seems difficult to prove.
\end{remark}

\begin{proof*}{\bf Proof of Theorem \ref{the2}}
The proof is very similar to that of Theorem \ref{the1}, this time we define $Z_n$ the same way as we have done in the previous proof, only now we define $\eta^Z_n(t)=B_n^{-1}Z_{nt}$ if $t=0,\frac{1}{n},\ldots 1$, and linear between these points, then \[\eta^Z_n(t)\Rightarrow W(t)\text{ for } t\in [0,1].\]
From Lemma \ref{rhoon12} we know that
$E_z(\rho_n)=O(\log n)$, and using similar arguments as in the previous proof we can show that
$B_n^{-1}\underset{0\le t\le 1}{\sup}||X_{[nt]}-Z_{[nt]}||$ tends to zero in probability as $n\to \infty$.  This proves the statement of our theorem.
\end{proof*}

\section{Proofs of lemmas}
The proofs of Theorem 1 and 2 are based on Lemma \ref{rhoon12}. In order to prove this, we need a few new definitions and lemmas.


\begin{definition}
For all $z\in K_N \cap Q$, consider a \rwwli such that $X_0=z$. Denote
\[\tau_z=\min \{ k\in \N | X_k \notin K_N \}\]
i.e. $\tau_z$ is the first exit time of the \rwwli from the set $K_N$.
\end{definition}

\begin{definition}\label{D:nun*}
Denote by $\nu_n$ the number of pure 1 blocks in the sequence
\[ \bs{1}_{X_0\in K_N}, \bs{1}_{X_1\in K_N}, \ldots, \bs{1}_{X_n\in K_N}\]
i.e. $\nu_n$ is the number of entrances into $K_N$ until time $n$.

We also define the slightly different random variable $\ol{\nu}_n$: it is the number of pure 1 blocks in the sequence
\begin{equation}\label{chisor}
\bs{1}_{x_0\in K_N}, \bs{1}_{x_1\in K_N}, \ldots, \bs{1}_{x_m\in K_N}
\end{equation}
where $m$ is chosen so that the number of 0s in the sequence (\ref{chisor}) equals $n$; i. e.
$\ol{\nu}_n$ is the number of entrances into $K_N$ until the first $n$ steps outside $K_N$.
\end{definition}
\begin{definition}\label{sbtb}
Denote by $B\subset Q$ a finite simply connected set, i.e., for all $x, y \in Q \setminus B$ there is a directed path from $x$ to $y$ in the subgraph on $Q \setminus B$. Let
\[S_B=\min \{k\in \N | X_k \in B\}\]
and
\[T_B=\min \{k\in \N | Y_k \in B\}.\]
These are called hitting times of $B$ by the processes $X_n$ and $Y_n$, respectively.
\end{definition}

\begin{lemma}\label{spratio}
Let $H\subset \l(Q\setminus K_N\r)$ be a finite set. Then there exists a constant $C_{H,K_N}>0$ (independent of $z$ and $n$), such that for all $z\in H$,
\[\frac{\ul{P}_z(T_{K_N}>n)}{\ul{P}_0(T_{\{0\}} > n)}>C_{H,K_N}\]
\end{lemma}

We shall denote by $E_z$ (and $\ul{E}_z$) expectations with respect to $P_z$ (and $\ul{P}_z$).


\begin{lemma}\label{pbsn}
If $d\ge 2$, then there is a $H\subset \l(Q\setminus K_N\r)$ finite set, and a $D_{H,K_N}> 0$ constant (independent of $z$ and $n$), such that for every $z\in Q$,
\[E_z(\nu_n)\le E_z(\ol{\nu}_n)\le \frac{D_{H,K_N}}{\underset{b\in H}{\min}\ul{P}_b(T_{K_N}>n)}\]
\end{lemma}

\begin{corollary}\label{nuon12}
If $d\ge 2$, then, as $n\to \infty$, $E_z(\ol{\nu}_n)=O(\log n)$.
\end{corollary}

\noindent Now we can start the proofs.

\begin{proof*}{\bf Proof of Lemma \ref{spratio}.}
We know that the following limit exists for all $z\in Q \setminus K_N$:
\[0 < \underset{n\to \infty}{\lim} \frac{P_z(T_{K_N}>n)}{\ul{P}_0(T_{\{0\}} > n)}<\infty\]
This is a special case of the Kesten-Spitzer Ratio Limit Theorem, for more details, see \cite{Spitzer}, page 165. From this, the lemma follows immediately.
\end{proof*}

%
%

\begin{proof*}{\bf Proof of Lemma \ref{pbsn}.}
By the definition it follows that $\ol{\nu}_n \ge \nu_n$, so $E_z(\ol{\nu}_n)\ge E_z(\nu_n)$, and we only need to prove the inequality for $\ol{\nu}_n$.

Let us denote by $S$ the series $X_0,X_1,X_2, \ldots, $, then create a new series $S'=\{\ol{X}_0,\ol{X}_1,\ldots \}$ by taking out every element of $S$ that are in $K_N$. Then we define $\phi: \Z^{+}\cup 0 \to \Z^{+}\cup 0$ such that the index of
$\ol{X}_i$ in the original series $S$ is $\phi(i)$ (it is easy to see that $\phi(i)\ge i$). Finally let the set $J$ be the following:
\[J=\{j: X_{\phi(j)+1}\in K_N\}\]
Then $J$ is the set of $j$ indexes such that we jump to $K_N$ from the element corresponding to $\ol{X}_j$ in $S$.
We define a complete system of events $A_i$, $i=0,1,...n$ and $A^c$:
\[A_i=\{i\in J, \{i+1,\ldots,n\} \text{ and } J \text{ are disjoint}\}\]
\[A^c=\{\{1,2,\ldots, n\}\text{ and } J \text{ are disjoint}\}\]
So the meaning of $A_i$ is that from the first $n$ steps outside $K_N$ in the series $S$, the last step when we jump to $K_N$ is the $i$th.
$A^c$ corresponds to the event that we do not get to $K_N$ in the first $n$ steps outside $K_N$ (thus $T_z(K_N)>n$).

We denote the complement of $K_N$ by $\ol{K_N}$. Then
\begin{eqnarray}\label{sumPAi1}
1&=&P_z(A^c)+\sum_{i=0}^n P_z(A_i)\ge \sum_{i=0}^n P_z(A_i)  \nonumber \\  &\ge& \sum_{i=0}^n \sum_{b\in \ol{K_N}}P_z(i\in J)P_z(\ol{X}_{i+1}=b |i \in J) P_b\l(S_{K_N}>n-i-1\r) 
\end{eqnarray}
Let $K_N^*$ be the set of those $x\in Q \cap K_N$ points, from which we can jump out of $K_N$ with positive probability. First we are going to deal with the case when there is such a $\tilde{b}\in (Q\setminus K_N)$ point that we can jump from every $x\in K_N^*$ to $\tilde{b}$ with positive probability (in 1 step). Then let $P_{m}=\underset{x\in K_N^*}{\min}P(x,\tilde{b})$, thus if we only take the term of $\tilde{b}$ in \eqref{sumPAi1}, in the summation to $b$, (and by using that $P_b\l(S_{K_N}>n-i-1\r)\ge P_b\l(S_{K_N}>n\r)$):
\begin{eqnarray}\label{sumPAi2}
1&\ge& \sum_{i=0}^n P_z(i\in J)P_{\tilde{b}}\l(S_{K_N}>n\r)\cdot P_z(\ol{X}_{i+1}=\tilde{b}| i\in J)\nonumber
\\&\ge& \sum_{i=0}^n P_z(i\in J)P_{\tilde{b}}\l(S_{K_N}>n\r) P_{m}
\end{eqnarray}
Moreover, $\sum_{i=0}^n P_z(i\in J)=E(\ol{\nu}_n)$, so by the choice $H=\{\tilde{b}\}$ and $D_{H,K_N}=1/P_m$ the statement of the lemma is true (we have also used that $P_b(S_{K_N}>n)$ is independent of the transition probabilities inside $K_N$, so $\ul{P}_b(T_{K_N}>n)=P_b(S_{K_N}>n)$).

The second case is when no such $\tilde{b}$ point exists. Now for all $x\in K_N^*$, we choose such a $b_x\in (Q\setminus K_N)$ that we can jump there with positive probability, i.e. $P(x,b_x)>0$. Let $H$ be the set of these $b_x$ points. It is evident from the construction that there is a constant $P_m>0$ such that jumping out of $K_N$, we get to a point in $H$ with greater or equal probability than $P_m$ ($P_m=\underset{x\in K_N^*}{\min}P(x,b_x)$ is a good choice). Then
\begin{align*}
1&\ge\sum_{i=0}^n \sum_{b\in \ol{K_N}}P_z(i\in J) P_z(\ol{X}_{i+1}=b |i \in J) P_b\l(S_{K_N}>n-i-1\r) \\
&\ge \sum_{i=0}^n P_z(i\in J)\sum_{b\in H}P_b\l(S_{K_N}>n\r) P_z(\ol{X}_{i+1}=b |i \in J)\\
&\ge \sum_{i=0}^n P_z(i\in J) \underset{b\in H}{\min} P_b\l(S_{K_N}>n\r) \sum_{b\in H} P_z(\ol{X}_{i+1}=b |i \in J)\\
&=\sum_{i=0}^n P_z(i\in J) \underset{b\in H}{\min} P_b\l(S_{K_N}>n\r) P_z(\ol{X}_{i+1}\in H|i \in J)\\
&\ge \sum_{i=0}^n P_z(i\in J) \underset{b\in H}{\min} P_b\l(S_{K_N}>n\r) P_m
\end{align*}
As in the previous case, we get the statement of our lemma with $D_{H,K_N}=1/P_m$.
\end{proof*}

\begin{proof*}{\bf Proof of Corollary \ref{nuon12}.}
The outline of the proof is similar to page 355-356., \cite{ED}, but here we prove it to random walk, instead of simple symmetric random walk. The idea to use this article comes from P\'{e}ter N\'andori.

First, we're going to suppose that the distribution has second moment. Then the probability of returning to the origin in the $n$th step is (by the Local Central Limit Theorem, \cite{Lawler}, page 25):
\begin{equation}\label{loktetel}
u(n)=1/\l(\l(2\pi\r)^{d/2} \sqrt{\det \Gamma}\cdot n^{d/2}\r)+o\l(1/n^{d/2}\r)=g/n^{d/2}+o\l(1/n^{d/2}\r)
\end{equation}
Here $\Gamma$ is the covariance-matrix (size $d \times d$), and $g=1/\l((2\pi)^{d/2} \sqrt{\det \Gamma}\r)$.

Denote by $R(n)$ the probability that the random walk does not return to the origin in $n$ steps, i.e. using Definition \ref{sbtb}, \[R(n)=\ul{P}_0(T_{\{0\}}>n),\] then
\[\sum_{k=0}^n u(k)R(n-k)=1\]
If $d\ge 3$, then
\[U=\sum_{n=0}^\infty u(n)<\infty\]
Let $R=\underset{n\to\infty}{\lim}R(n)$, which exists because $R(n)$ is bounded and monotone. For $1\le k\le n$
\[R(n-k)\sum_{i=0}^k u(i)+\sum_{i=k+1}^n u(n)\ge 1\]
If $k\to \infty$ and $n-k\to \infty$, then
\[R(n-k)\cdot U\ge 1+o(1)\]
so \[R\ge \frac{1}{U}\]
But $R(n)$ is monotone decreasing, so \[R(n)\ge \frac{1}{U}.\]
If $d=2$, then by substituting \eqref{loktetel} in the place of $u(k)$, we get
\[u(0)+u(1)+\ldots u(n)=g\log n(1+o(1))\]
using that $R(n)$ is  decreasing
\begin{equation}\label{logle}
R(n)g\log n\le 1+o(1)
\end{equation}
We know that for $0<k<n$,
\[R(n-k)[u(0)+\ldots u(k)]+u(k+1)+\ldots +u(n)\ge 1\]
If $k$ and $n$ tends to infinity, then
\[R(n-k)\cdot g\log(k) (1+o(1))+g(1+o(1))\log \frac{n}{k}\ge 1\]
Let $k=n-[n/\log n]$, then
\[R(n-k)g\log(n-k)(1+o(1))+o(1)\ge 1\]
This and \eqref{logle} yields
\[R(n)=\frac{1+o(1)}{g\log n}\]

If the distribution is B-type, then in $d$ dimensions, according to the Local Central Limit Theorem in \cite{Rv} (see Appendix), there is such a $B_n$ series and $c=g(0)$ positive constant such that
\[u(n)=\frac{c}{B_n^{d}}+o\l(\frac{1}{B_n^{d}}\r)\]
We see that $\underset{n\to \infty}{\lim}\frac{B_n}{\sqrt{n}}>0$. Let us call $C_n=\sum_{k=0}^n u(k)=\sum_{k=0}^n \frac{c}{B_n^2}$, then it is monotone increasing. For $d=2$, we have $C_n=O(\log n)$, and for $d\ge 3$, $C_n=O(1)$.

Now, depending on whether $\underset{n\to \infty}{\lim}C_n=U<+\infty$ or $\underset{n\to \infty}{\lim}C_n=+\infty$, we have two cases. In the first case, using the same argument as before, we can show that $R(n)\ge \frac{1}{U}$.
In the second one, we also use the same argument as previously: using the fact that $\sum_{k=0}^n u(k)R(n-k)=1$, we can easily see that
\begin{equation}\label{Cnle}
R(n)C_n\le 1+o(1)
\end{equation}

We know that for $0<k<n$,
\[R(n-k)[u(0)+\ldots u(k)]+u(k+1)+\ldots +u(n)\ge 1\]
If $k$ and $n$ tends to infinity, then
\[R(n-k)C_k(1+o(1))+(C_n-C_k)\ge 1\]
Let $k=n-[n/C_n]$, then \[C_n-C_k\le c(\log n-\log(n-[n/C_n]))(1+o(1))=c\log\l(\frac{1}{1-1/C_n}\r)(1+o(1))=o(1).\]
On the other hand, we see that \begin{align*}
C_{k}-C_{n-k}&\le (\log(n-[n/C_n])-\log([n/C_n]))(1+o(1))=\\ &=\log(C_n-1)(1+o(1))=o(C_{n})=o(C_{k})=o(C_{n-k}),
\end{align*}
therefore we can write
\[R(n-k)C_{n-k}(1+o(1))+o(1)\ge 1\]
This and \eqref{logle} yields
\[R(n)=\ul{P}_0(T_{\{0\}}>n)=\frac{1+o(1)}{C_n}\]

From these, Lemma \ref{spratio}. and Lemma \ref{pbsn}., we get that
\[\underset{n\to \infty}{\ol{\lim}} \ul{P}_0(T_{\{0\}}>n)E_z(\ol{\nu}_n)\le \underset{n\to \infty}{\ol{\lim}} \frac{D_{H,K_N}\ul{P}_0(T_{\{0\}}>n)}{\underset{b\in H}{\min}\ul{P}_b(T_{K_N}>n)}<\infty \]
\[E_z(\ol{\nu}_n)=O(\log n)\]

\end{proof*}

\begin{proof*}{\bf Proof of Lemma \ref{rhoon12}.}
\begin{definition}
Let us denote the time spent outside of $K_N$ between the $i-1$th and $i$th visit to $K_N$ by $\eta_{i-1}$. This is the same as the length of the $i-1$th zero block in the series $\bs{1}_{x_0\in K_N}, \bs{1}_{x_1\in K_N}, \ldots $. Let $\xi_{i-1}$ be the time spent in $K_N$ by $X_N$ during the $i$th visit, which equals the length of the $i$th ones block in the previous series.
\end{definition}
The difficulty in the proof is that $\xi_i$ are dependent, so they are hard to handle directly. For this reason, we're going to work with the conditional probabilities in function of the entrance/exit points.
\begin{definition}
Let $u_k$ be the point of the $k$th entrance to $K_N$ ($u_k\in K_N$), and $v_k$  the last point in $K_N$ before the $k$th exit from $K_N$. Let us denote $U_n=\{u_k\}_{k=1}^n$ and $V_n=\{v_k\}_{k=1}^n$. Then $U_n$ and $V_n$ are random vectors taking values in $K_N^n$.
\end{definition}
There exists such a $K<\infty$ that
$E_z(\xi_k | u_k, v_k)\le K$
because the expected value only depends on the value of $u_k$ and $v_k$ (but does not depends on $k$), and there are at most $\l|K_N\r|^2$ $u,v$ pairs, and the expected value is bounded for all of them, so we can choose $K$ as the maximum of these. In the following, we'll use the shorthand notation for summation $\sum_{U_n,V_n}  \overset{def}{=}\sum_{U_n\in K_N^n,V_n\in K_N^n}$,
\begin{align*}
E_z(\rho_n)&\le E_z(\sum_{i=1}^{\ol{\nu}_n}\xi_i)=\sum_{k=0}^n E(\sum_{i=1}^{k}\xi_i| \ol{\nu}_n=k)P(\ol{\nu}_n=k)=\\
&=\sum_{k=0}^n \sum_{i=1}^{k} \sum_{U_n,V_n} E(\xi_i |\ol{\nu}_n=k, u_i(U_n),v_i(U_n))P(\ol{\nu}_n=k,U_n,V_n)\le\\
&\le \sum_{k=0}^n \sum_{U_n,V_n} K\cdot k \cdot P(U_n,V_n,\ol{\nu}_n=k)=K  \sum_{k=0}^n k \cdot P(\ol{\nu}_n=k)=K\cdot E(\ol{\nu}_n)
\end{align*}

Here we have used that $E(\xi_i | u_i,v_i, \ol{\nu}_n=k)=E(\xi_i | u_i,v_i)\le K$, this is true because $\ol{\nu}_n$ is independent of the time spent inside $K_N$. Now it is clear that $E_z(\rho_n)=O(\log n)$.
\end{proof*}

\section{Remarks}
\begin{enumerate}
\item
This paper is based on the work done for D. Paulin's bachelor's thesis in 2009 at Budapest University of Technology and Economics.
\item
The conditions for our Theorem 1 are quite general. With some technical work, it could be easily shown that Theorem 1 also holds for periodic random walks, the same way it is done in \cite{SzT}. 
\item
In relation to Theorem 2, we conjecture that, analogously to the finite horizon case, the weak limit of the  locally perturbed Lorentz-process with infinite horizon is the same  Brownian motion as it is that of the periodic one.
\item
If for the distribution of the jumps in the random walk, \[\underset{k\in \R^+}\sup (k. \text{ moment exists})=\alpha,\] with $\alpha\in [1,2)$, and the tail of the random walk satisfies some criteria (see \cite{Whitt} page 114 Theorem 4.5.1),
then we will converge to a stable distribution with parameter $\alpha$, and the normalizing factor will not be $n^{1/2}$, but $n^{1/{\alpha}}L_0(n)$, $L_0$ being a slowly varying function.
We can apply the local limit theorem for stable distributions, Theorem 6.1 of \cite{Rv}, to the returns to the origin,
$P(Y_n^{(1)}=0)=O(n^{-1/{\alpha}})$.
Using the same theorem we get that in two dimensions, if the jump belongs to the domain of attraction of the 2 dimensional stable law with parameter $\alpha$, then $P\l(Y_n^{(2)}=0\r)= O(n^{-2/{\alpha}})$. $\alpha\in [1,2)$, so $\sum_{n=0}^{\infty}P\l(Y_n^{(2)}=0\r)<\infty$, and thus the random walk is transient. From this, the convergence to a L\'{e}vy - process can be proven.
\end{enumerate}

\begin{acknowledgements}
The authors are grateful to P\'eter N\'andori for his constant support and invaluable comments. We also express our sincere thanks to the referees for their careful reading of the manuscript and their valuable remarks.
\end{acknowledgements}

\section{Appendix}
In what follows we summarize some notions and theorems related to multidimensional limit theory of sums of iid random vectors, in particular some global and local limit theorems and domains of attractions. These are typically not included in textbooks on probability, and moreover, the pioneering work of Rvaceva on local theorems is not easily available.

\begin{definition}\label{D:Dattr}
Let $\{\xi(n)\}_n$ be a sequence of iid random vectors, with distribution function $F(x)$. Then if there are suitably chosen constants $C(n)>0$, real vectors $d(n)$ such that $s_n=\sum_{k=1}^n\xi_k/C(n)-d(n)$ converges in distribution to a non-degenerate probability distribution $R(x)$, then
\begin{enumerate}
\item
$R$ is called a stable law;
\item and we say that $F(x)$ belongs to the domain of attraction of $R(x)$.
\end{enumerate}
\end{definition}

Our main interest in this paper is the domain of attraction of the gaussian law, therefore below we also restrict our attention to it.

Theorem 4.1 of \cite{Rv} describes the domain of attraction of the normal law for random vectors (below ' means matrix transpose):
\begin{quote}
\textbf{Theorem 4.1}(\cite{Rv}). $P$ belongs to the domain of attraction of the non-degenerate normal law with characteristic function $\exp(-Q(t)/2)$ if and only if:
\begin{align*}
&(1)\hspace{10mm} R^2\int_{|x|>R}\dd P(x)/\int_{|x|<R}x'x\dd P(x) \to 0 \text{ as } R\to \infty,\\
&(2)\hspace{10mm} \int_{|x|<R}(t'x)^2 \dd P(x)/\int_{|x|<R}(u'x)^2\dd P(x)\to Q(t)/Q(u) \text{ as } R\to \infty,
\end{align*}
for arbitrary $t,u \in \mathbb R^d$.
\end{quote}

Further if the jump distribution $\ul{P}$ of the random walk is B-type, then we define $\eta^{Y'}_n(t)=B_n^{-1}Y_{nt}$ for $t=0,\frac{1}{n},\frac{2}{n},\ldots,1$ and take its piecewise linear extension. It is well-known that $\eta^{Y'}_n\Rightarrow W_{\sum}(t)$, a $d$-dimensional Wiener process (in fact, Skorohod proved that when the CLT holds, then this functional central limit theorem holds too; for more details, see \cite{Sko} and \cite{Whitt}, page 115-118). The covariance matrix $\sum$ is determined by Theorem 4.1 of \cite{Rv}.

\begin{definition}\label{D:Ltype}
As a special case, we call a B-type random variable (or random vector) L-type if $B_n=\sqrt{cn\log n}$. (N. B.: this scaling is used in the weak limit of the planar, infinite-horizon Lorentz-process; cf. \cite{B} and \cite{SzV2007}). Now $L(x)\sim 2c\log x$. If we denote by $S_n^*$ the  sum of $n$ iid one dimensional L-type variables, then $\frac{S_n^*}{\sqrt{cn\log n}} \overset{d}{\to} N(0,1)$.
\end{definition}

Finally we recall Rvaceva's Local Limit Theorem, Theorem 6.1 of \cite{Rv}:

Let $\{\xi(n)\}_n$ be a sequence of iid $\mathbb Z^p$-valued random vectors and let $P(x)=Pr(\xi(n)=x)$. Let $P(n;z)=Pr(s(n)=z)$, where $s(n)=\xi(1)+\ldots+\xi(n)$, and $g(x)$ be the density of a certain stable distribution $G$.

\textbf{Theorem 6.1 of \cite{Rv}} In order that for some suitably chosen constant vectors $a(n)$ and positive constants $B(n)$ the relation
\[R(n)=B^p(n) P(n;z)-g[[z-a(n)]/B(n)]\to 0\]
hold uniformly with respect to $z$, it is necessary and sufficient that the distribution of $\xi(n)$
\begin{enumerate}
\item[1)] belong to the domain of attraction of $G$, and
\item[2)] be a 1 - lattice distribution
\item[] The second condition is equivalent to each of the following:
\item[2')]The greatest common divisor of the volumes of $p$ - dimensional simplexes a $p+1$ vertices of which lie at points with $P(x)>0$ is $1/p!$
\item[2'')]the lattice generated by all vectors $(x-y)$ such that $P(x)>0<P(y)$ coincides with the lattice of all integral points of the $p$ - dimensional space
\end{enumerate}
\begin{remark}\label{R:degen}
The series $B(n)$ in this theorem may be different from $C(n)$ in Definition \ref{D:Dattr}. If a two dimensional RW jumps independently along the axes, and it has a distribution with finite second moment along one axis and an L-type distribution along the other, then $C(n)\sim \sqrt{ n\log n}$ and $B(n)\sim \sqrt{n}\log^{1/4} n$.
\end{remark}


\begin{thebibliography}{99}
\bibitem{B} Bleher, P. M.: {Statistical Properties of Two-Dimensional
  Periodic Lorentz Gas with Infinite Horizon}. J. of Stat.\ Physics,
  66(1):315--373, 1992.

\bibitem{ChD}
Chernov, N.  and Dolgopyat, D.:
 {Hyperbolic billiards and statistical physics}
Proc. of International Congress of Mathematicians (Madrid, Spain, August 2006), Vol. II, Euro. Math. Soc., Zurich, 2006, pp. 1679-1704.

\bibitem{ChD2009} Chernov, N., and Dolgopyat, D.: {Anomalous current in periodic Lorentz gases with infinite horizon},
Uspekhi Mat. Nauk, 64:4 (2009), 73-124 (in Russian),
Russ. Math. Surveys, 64 (2009) 651-699

\bibitem{ChD2010} Chernov, N., and Dolgopyat, D.:
{Lorentz gas with thermostatted walls}, pp. 50, submitted

\bibitem{DSzV}
Dolgopyat, D., Sz\'asz, D. and Varj\'u, T.: {Limit Theorems for Locally Perturbed Planar Lorentz Processes},
Duke Math. Journal, 148: 459-499, 2009

\bibitem{ED}
Dvoretzky, A.  and Erd\H{o}s, P.:
{Some problems on random walk in space. }
Proc. $2^{\text{nd}}$ Berkeley Sympos. Math. Statis. Probab., pp. 353-367 (1951)


\bibitem{Gnedenko}
Gnedenko, B.V.  and Kolmogorov, A.N.:
 {Limit distributions for sums of independent random variables. }
Revised Edition (1968), Translated by K.L.Chung. Afterword: J. L. Doob. Addison-Wesley series in statistics. ASIN: B000WSX412

\bibitem{Skew}
Harrison, J. M.  and Shepp, L. A.
{On Skew Brownian Motion}
Ann. Probab. Volume 9, Number 2 (1981), 309-313.

\bibitem{Independent}
Ibragimov, I. A. and Linnik, Yu. V.:
 {Independent and Stationary Sequences of Random Variables}
1971 Wolters-Noordhoff Publishing Groningen

\bibitem{Lawler}
Lawler, G. F.  and Limic, V.:  {Random walk: a modern introduction.} To be published by Cambridge University Press. Electronic version available at http://www.math.uchicago.edu/$\sim$lawler/srwbook.pdf.

\bibitem{Lindvall}
Lindvall, T.: {Weak convergence of probability measures and random functions in the function space $D[0, \infty)$}
J. Apl. Probab. 10 (1973), 109-121


\bibitem{MS08}
Marklof, J. and Str\"ombergsson, A.: {Kinetic transport in the two-dimensional periodic Lorentz gas}, Nonlinearity 21 (2008), 1413-1422.

\bibitem{Rv}
Rvaceva, E.
{On the domains of attraction of multidimensional distributions,}
Selected Transl. Math. Stat. Prob., 2 (1962), 183-207.

\bibitem{Sko}
{
Skorokhod, A. V.:
{Limit Theorems for Stochastic Processes}
Theory Probab. Appl. 1, 261 (1956)
}

\bibitem{Spitzer}
Spitzer, F.: {Principles of Random Walk}, 2nd edition, Sringer-Verlag,1976,
ISBN-10: 0387951547,
ISBN-13: 978-0387951546

\bibitem{Stone}
Stone, C.
{A Local Limit Theorem for Nonlattice Multi-Dimensional Distribution Functions,}
Annals of Math. Statistics,  36 (1965), 546-551

\bibitem{SzT}
Sz\'asz, D. and Telcs, A.:
 {Random Walk in an Inhomogeneous Medium with Local Impurities. }
Journal of Statistical Physics, Vol. 26, No. 3, 1981\\
ISSN: 0022-4715 (Print) 1572-9613 (Online)

\bibitem{SzV2004}
Sz\'asz, D. and Varj\'u,T.:
 {Local Limit Theorem and Recurrence for the Planar Lorentz Process}
Ergodic Theory and Dynamical Systems (2004), 24, 257-278

\bibitem{SzV2007}
Sz\'asz, D. and Varj\'u,T.:
 {Limit Laws and Recurrence for the Planar Lorentz Process with Infinite Horizon}
Journal of Statistical Physics (2007) 129: 59-80

\bibitem{Whitt}
Whitt, W.:
{Stochastic-Process Limits,}
Springer Series in Operations Research, Springer-Verlag, New York, 2002.

\end{thebibliography}
\end{document}